\documentclass{article}
\usepackage[utf8]{inputenc}
\usepackage{mathtools}
\usepackage{amsmath,amssymb,amsbsy}

\newcommand{\zet}{\mathbb{Z}}

\newcommand{\qed}{\hfill \rule{.1in}{.1in}}

\newtheorem{thm}{Theorem}[section]
\newtheorem{lem}[thm]{Lemma}

\newtheorem{cor}[thm]{Corollary}
\newtheorem{prob}[thm]{Problem}
\newtheorem{dfn}[thm]{Definition}

\newtheorem{conj}[thm]{Conjecture}

\title{Disjoint zero-sum subsets in Abelian groups and its application - survey}

\author{Sylwia Cichacz
\normalsize \\AGH University, \vspace{2mm} Poland\\
}

\begin{document}
\maketitle
\begin{abstract}
We provide a summary of research on disjoint zero-sum subsets in finite Abelian groups, which is a branch of additive group theory and combinatorial number theory.

An orthomorphism of a group $\Gamma$ is defined
as a bijection   $\varphi$ $\Gamma$ such that the mapping  $g \mapsto  g^{-1}\varphi(g)$ is
also bijective. In 1981, Friedlander, Gordon, and Tannenbaum conjectured that when $\Gamma$ is Abelian, for any $k \geq 2$ dividing $|\Gamma| -1$,
there exists an orthomorphism of $\Gamma$ fixing the identity and permuting the remaining elements as products of disjoint $k$-cycles. Using the idea of disjoint-zero sum subset we provide a solution of this conjecture for $k=3$ and $|\Gamma|\cong 4\pmod{24}$.

We also present some applications of zero-sum sets  in graph labeling. 
\end{abstract}

 \section{Introduction}
Assume $\Gamma$ is a finite Abelian group with the operation denoted by~$+$.
  For convenience
we will write $ka$ to denote $a + a + \ldots + a$ where the element $a$ appears $k$ times, $-a$ to denote the inverse of $a$, and
we will use $a - b$ instead of $a+(-b)$.  Moreover, the notation $\sum_{a\in S}{a}$ will be used as a short form for $a_1+a_2+a_3+\dots$, where $a_1, a_2, a_3, \dots$ are all elements of the set $S$. The identity element of $\Gamma$ will be denoted by $0$. A subset $S$ of $\Gamma$ is referred to as a \textit{zero-sum subset} if $\sum_{g\in S}g=0$. Zero-sum problems typically examine the conditions under which given sequences contain non-empty zero-sum subsequences with particular properties \cite{Gao}. The concept was introduced by Erd\H{o}s, Ginzburg, and Ziv, who demonstrated that $2n-1$ is the smallest integer such that every sequence over a cyclic group 
$\mathbb{Z}_n$   has a zero-sum subsequence of length 
$n$ \cite{ErdGinZiv}. This finding sparked extensive research in the area.

For example, a well-studied aspect of zero-sum problems for graphs is zero-sum Ramsey theory, which can be framed as follows: What is the smallest number $n$ such that any complete graph $K_m$ ($m\geq n$) with edges labeled by elements of a finite group 
$\Gamma$ contains a subgraph of a specified type where the total weight of the edges is zero in $\Gamma$?

Another problem for zero-sum sets is the issue of disjoint subsets in $\Gamma$. This approach was inspired by research on Steiner triples and initiated by Skolem in 1957, who posed the following question \cite{Skolem57}:

\textit{
Is there a partition of the set of nonzero elements of the cyclic group $\mathbb{Z}_n$   into three-element subsets such that the sum of the elements in each subset is $0\pmod n$, for $n\equiv 1\pmod 6$?
}

This question received an affirmative answer \cite{ref_Han, Skolem57, Skolem58} and served as the starting point for research of Steiner triples in Abelian groups \cite{Tannenbaum3} and a more general problem posed by Tannenbaum:
\begin{prob}[\cite{Tannenbaum1}]\label{problemT}
Let $\Gamma$ be an abelian group of order $m$. Let 
$t$ and $m_i$ for $i\in\{1,\ldots,t\}$ be positive integers such that $\sum_{i=1}^t m_i = m - 1$. Let  $w_i\in\Gamma$ for $i \in \{1,\ldots t\}$, Determine when there exists a partition of the set  $\Gamma^*=\Gamma\setminus\{0\}$ subsets $\{S_i\}_{i=1}^t$   such that $|S_i| = m_i$ and $\sum_{s\in S_i}s = w_i$   for each $i\in\{1,\ldots,t\}$. 
\end{prob}
If a respective subset partition of $\Gamma^*$ exists, then we say that $\{m_i\}_{i=1}^t$ is \textit{realizable} in $\Gamma^*$ with $\{w_i\}_{i=1}^{t}$. In such a case, we say that $\{S_i\}_{i=1}^t$ \textit{realizes} $\{m_i\}_{i=1}^t$ in $\Gamma^*$ with $\{w_i\}_{i=1}^{t}$. 

The most extensively studied case in the literature is the \textit{Zero-Sum Partition} (ZSP) problem, which occurs when 
$w_i=0$ for every  {$i\in\{1,\ldots,t\}$}. In this paper, we will present a survey on this concept as well as an application of it to graph labeling problems. Note that many results on disjoint zero-sum sets have been obtained independently due to their application to various graph problems.

  \section{Preliminaries}
For positive integers $a$ and $b$ such that $a<b$ let $[a,b]=\{a,a+1,\ldots,b\}$.

 Recall that any group element $\iota\in\Gamma$ of order 2 (i.e., $\iota\neq 0$ and $2\iota=0$) is called an \emph{involution}. Let us denote the set of involutions in $\Gamma$ by $I(\Gamma)$. The sum of all elements of a group $\Gamma$ is equal to the sum of its involutions and the identity element.
Therefore $\sum_{g\in \Gamma}g= \iota$ if $|I(\Gamma)|=1$
 (where $\iota$ is the involution) and $\sum_{g\in \Gamma}g= 0$, otherwise (see \cite{ref_ComNelPal}, Lemma 8).

\section{Zero-sum partition of Abelian group}
The case most studied in the literature is the {\em Zero-Sum Partition (ZSP)} problem, i.e., when $w_i=0$ for every  {$i\in[1,t]$} in Problem~\ref{problemT}. Note that,
 for an Abelian group $\Gamma$ with $|I(\Gamma)|=1$, $\Gamma^*$ does not admit zero-sum partitions. Moreover, given a group $\Gamma$, if $\{m_i\}_{i=1}^t$ is realizable in $\Gamma^*$ (with $w_i=0$ for every {$i\in[1,t]$}), then necessarily $m_i\geq 2$ for every {$i\in[1,t]$}. It was proved that this condition is sufficient when:
 the group $\Gamma$ is cyclic of odd order  \cite{FrieGorTan,KLR}, for groups $(\zet_p)^n$, where $p>2$ is a prime \cite{Fukuchi},  for any Abelian group of odd order  \cite{Tannenbaum1,Zeng}, and finally for $\Gamma$ having exactly three involutions \cite{Zeng}.
 
  Let us generalize this situation with the following definition.

\begin{dfn}[\cite{CicrSuch}]Let $\Gamma$ be a finite Abelian group of order $m$. We say that $\Gamma$ has {\em $x$-Zero-Sum Partition Property ($x$-ZSPP)} if, for every positive integer $t$ and every integer partition $\{m_i\}_{i=1}^t$ of $m-1$ (i.e., $\sum_{i=1}^t m_i=m-1$) with $m_i \geq x$ for every  {$i\in[1,t]$}, there exists a subset partition $\{S_i\}_{i=1}^t$ of $\Gamma^*$ (i.e., subsets $S_i$ are pairwise disjoint and their union is $ \Gamma^*$) with $|S_i| = m_i$ and $\sum_{s\in S_i}s = 0$ for every {$i\in[1,t]$}. 
\end{dfn}

The following theorem was first conjectured by Kaplan,  Lev and Roditty in \cite{KLR} who also demonstrated its necessity, and was later proven by Zeng in \cite{Zeng}.
\begin{thm}[\cite{Zeng}]\label{Zeng}Let $\Gamma$ be a finite Abelian group of order $m$. $\Gamma$ has $2$-ZSPP if and only if $|I(\Gamma)|\in\{0,3\}$. 
\end{thm}

It was shown independently by a few authors that 
$\Gamma\cong (\zet_2)^n$, $n>1$ has $3$-ZSPP.

\begin{thm}[\cite{CaccJia,Egawa,Tannenbaum2}]\label{Sylow}Let $\Gamma\cong (\zet_2)^n$ for some integer $n$, with $n>1$, then $\Gamma$ has $3$-ZSPP. 
\end{thm}

The above theorems confirm, for the case of $|I(\Gamma)|=3$ or $\Gamma\cong (\zet_2)^n$, $n>1$ (Note that in this case $I(\Gamma)=\Gamma^*$, so $R =\Gamma \setminus (I(\Gamma)\cup\{0\})$), the following conjecture stated by Tannenbaum. 
\begin{conj}[\cite{Tannenbaum1}]\label{conjectureT}
Let $\Gamma$ be a finite Abelian group of order $m$ with $|I(\Gamma)|>1$. Let $R =\Gamma \setminus (I(\Gamma)\cup\{0\})$. For every positive integer $t$ and every integer partition $\{m_i\}_{i=1}^t$ of $m-1$, with $m_i \geq 2$ for every {$i\in [1, |R|/2]$}, and {$m_i \geq 3$} for every {$i\in[|R|/2+1, t]$}, there is a subset partition $\{S_i\}_{i=1}^t$ of $\Gamma^*$ such that $|S_i| = m_i$ and $\sum_{s\in S_i}s = 0$ for every {$i\in[ 1, t]$}. 
\end{conj}
Cichacz and Suchan showed recently (\cite{CicrSuch2}) that Conjecture \ref{conjectureT} is false in general, and it is only true in the cases covered by Theorems \ref{Zeng} and \ref{Sylow}.

\begin{thm}[\cite{CicrSuch2}]\label{Skol}
Let $\Gamma$ be a finite Abelian group with $|I(\Gamma)|>1$. Let $R =\Gamma^* \setminus I(\Gamma)$. For any positive integer $t$ and an integer partition $\{m_i\}_{i=1}^t$ of $|\Gamma^*|$, with {$m_i \geq 2$ for all $i$, $i\in[1, |R|/2]$, and {$m_i \geq 3$} for all $i$, $i\in[|R|/2+1,t]$}, there is a subset partition $\{S_i\}_{i=1}^t$ of $\Gamma^*$ such that $|S_i| = m_i$ and $\sum_{s\in S_i}s = 0$ for all  $i\in[1,t]$ if and only if $|I(\Gamma)|\in\{3,|\Gamma^*|\}$.\end{thm}

For every Abelian group $\Gamma$ with more than one involution and large enough order, first Cichacz and Tuza~\cite{CicTuz} showed that it has $4$-ZSPP. This result was improved by Cichacz and Suchan~\cite{CicrSuch}, for every $\Gamma$ of order $2^{n}$ for some integer $n>1$ such that $|I(\Gamma)|\neq 1$, with the following theorem. 
\begin{thm}[\cite{CicrSuch}]\label{SK}Let $\Gamma$ be such that $|I(\Gamma)|>1$ and $|\Gamma|=2^n$ for some integer $n>1$. Then $\Gamma$ has {$3$-ZSPP.}
\end{thm} 
Later M\"{u}yesser and Pokrovskiy showed that if $\Gamma$ has more than one involution and large enough order, then it indeed has $3$-ZSPP~\cite{MP}. Finally Cichacz and Suchan proved the following:
\begin{thm}[\cite{CicrSuch2}]\label{SK2}Let $\Gamma$ be such that $|I(\Gamma)|>1$. Then $\Gamma$ has {$4$-ZSPP.}
\end{thm}

Therefore Conjecture~\ref{conjectureT} posed by Cichacz \cite{CicZ} is still open.
\begin{conj}[\cite{CicZ}]\label{conjecture}Let $\Gamma$ be a finite Abelian group with $|I(\Gamma)|>1$. Then $\Gamma$ has $3$-ZSPP. 
\end{conj}

However, it has achieved important progress towards proving Conjecture~\ref{conjecture}. Recall that any group $\Gamma$ can be factorized as $\Gamma\cong L\times H$, where $L$ is the Sylow $2$-group of $\Gamma$ and the order of $H$ is odd. In this context, towards proving Conjecture~\ref{conjecture}, Cichacz and Suchan left open only the question if $\Gamma$ has not only $4$-ZSPP, but also $3$-ZSPP, in the case where $(|H| \bmod{6}) = 5$, for the other cases (i.e $(|H| \bmod{6})\in\{1,3\}$) the group $\Gamma$ has $3$-ZSPP  \cite{CicrSuch2}. 

Note that the motivation for Conjecture~\ref{conjecture} was the following result:
\begin{thm}[\cite{CicZ}]Let $\Gamma$ be a finite Abelian group with $|I(\Gamma)|>1$. If $m>2$ divides $|\Gamma|$, then there is a partition of $\Gamma$ into pairwise disjoint subsets $S_1,S_2,\ldots,$ $S_t$, such that  with $|S_i| = m$ and $\sum_{s\in S_i}s = 0$ for every {$i\in[1,t]$}.
\end{thm}

\subsection{Skolem partition of Abelian groups}
A Skolem partition is a concept related to combinatorial design, particularly involving sequences and partitions of sets with certain properties. It is named after the Norwegian mathematician Thoralf Skolem \cite{Skolem57}. 

A \textit{Skolem sequence} of order $n$ is a sequence of $2n$ elements that consists of the numbers 
$1,1,2,2,\ldots,n,n$ arranged in such a way that for each $k$ (where 
$1\leq k \leq n$), the two occurrences of $k$ are exactly $k$ positions apart. For example,  the sequence $42324311$ is a Skolem sequence.

Skolem provides the method to derive Steiner triple system from a Skolem sequence \cite{Skolem58} and acted as a foundation for a concept introduced by Tannenbaum. Namely, we call a $6$-subset $C$ of an Abelian group $\Gamma$ \textit{good} if $C = \{c, d,-c -d,-c,-d, c + d\}$ for some $c$ and $d$ in $\Gamma$. Notice that the sum of elements of a good $6$-subset is $0$. Moreover, it can be partitioned into three zero-sum $2$-subsets or two zero-sum $3$-subsets.


 The following definition was given by Tannenbaum~\cite{Tannenbaum1}.

\begin{dfn}\label{dfn:skolem}Let $\Gamma$ be a finite Abelian group of order $m=6k+s$ for a non-negative integer $k$ and $s\in\{1,3,5\}$. A partition of $\Gamma^*$ into $k$ good $6$-subsets and $(s-1)/2$ zero-sum $2$-subsets is called a {\em Skolem partition of $\Gamma^*$}.
\end{dfn}

Note that, if $\Gamma$ is an Abelian group of order $m$ such that every integer partition $\{m_i\}_{i=1}^t$ of $m-1$ with $m_i \in \{2, 3\}$ for every $i\in[1,t]$ is realizable in $\Gamma^*$, then also every integer partition $\{m'_i\}_{i=1}^{t'}$ of $m-1$ with $m'_i \geq 2$ for every $i\in[1,t']$ is realizable. Similarly, it is easy to see that if a Skolem partition of $\Gamma^*$ exists, then every integer partition $\{m_i\}_{i=1}^t$ of $m-1$ with $m_i \in \{2, 3\}$ for every $i\in[1,t]$ is realizable. So the following theorem by Tannenbaun indeed offers a stronger version of $2$-ZSPP.

\begin{thm}[\cite{Tannenbaum1}]\label{Tannenbaum1}Let $\Gamma$ be a finite Abelian group such that $|I(\Gamma)|=0$, then $\Gamma^*$ has a Skolem partition. 
\end{thm}

It is worth mentioning that 6-good subsets have significant applications in proving the $x$-ZSPP \cite{CicrSuch,CicrSuch2,CicTuz,Zeng} of groups, thus 
in \cite{CicrSuch2} the authors generalized  Definition \ref{dfn:skolem}. 

Let $S$ be a subset of cardinality $m=6k+s$, for a positive integer $k$ and $s\in\{0,2,4\}$, of a finite Abelian group $\Gamma$. A partition of $S$ into $k$ good $6$-subsets and $s/2$ zero-sum $2$-subsets is called a {\em Skolem partition of $S$}. It is known that $\zet_m^*\setminus I(\zet_m)${, for $m=6k+s$ with any positive integer $k$,} has a Skolem partition when $s\in\{0,4\}$, but the situation is slightly different for $s=2$. For $m\equiv 2$ or $8\pmod{24}$, there exists a Skolem partition, whereas for $m\equiv 14$ or $20\pmod{24}$, such a partition does not exist (see Lemma 4 in \cite{Tannenbaum1}).

Note that there exist groups $\Gamma$ with $|I(\Gamma)|>1$ such that the set $R =\Gamma^* \setminus I(\Gamma)$ has a Skolem partition. For instance, let $\Gamma\cong(\zet_2)^{\eta}\times H$  for some natural number $\eta>1$ and $|H|\equiv 1\pmod 6$ \cite{CicrSuch2}.

We finish this section with the open problem.
\begin{prob}[\cite{CicrSuch2}]\label{conjectureSK}Characterize finite Abelian groups $\Gamma$ for which $R =\Gamma^* \setminus I(\Gamma)$ has a Skolem partition.
\end{prob}

\section{Disjoint zero-sum subsets}
Recall that for an Abelian group $\Gamma$ with $|I(\Gamma)|=1$, $\Gamma^*$ does not admit zero-sum partitions. Therefore for such groups one can consider the set $R=\Gamma^*\setminus I(\Gamma)$.

For  cyclic groups it was proven:
\begin{thm}[\cite{Tannenbaum1,Zeng}]\label{Zeng2}
{Let} $n$ be even positive natural number and  $m,l$ be natural numbers such that $3m+2l = n-2$. Then the set $R=\zet_n\setminus\{0,\frac{n}{2}\}$ can be partitioned into $m$ zero-sum $3$-sets $A_1,A_2,\ldots,A_{m}$ and $l$ zero-sum $2$-sets $B_1,B_2,\ldots,B_l$.  
\end{thm}

For general groups with one involution a weaker result was obtained recently:
\begin{thm}[\cite{Cicir}]   \label{zero-sum}

Let $\Gamma$ be of order\/ $n$, $|I(\Gamma)|=1$, $\iota_{\Gamma}$ be the involution in $\Gamma$ and consider any integers\/
 $r_1,r_2,\dots,r_t$ with\/ $n -2=r_1 + r_2 + \ldots + r_t$
  and {$r_i \geq 4$} for all\/ $i\in[1,t]$.
Then
 there exist pairwise disjoint zero-sum subsets\/ $S_1, S_2,\ldots , S_t$ in\/ $\Gamma\setminus\{0,\iota_{\Gamma}\}$
 such that\/ $|S_i| = r_i$ for all\/ $i\in[1,t]$.
\end{thm}
 For $\Gamma$ with $|I(\Gamma)|=1$, we believe that Theorem~\ref{Zeng2} holds not only for cyclic groups and we state the following conjecture.


\begin{conj}[\cite{Cicir}]\label{moja}
Let $\Gamma$ be a finite Abelian group of order $m$ with $|I(\Gamma)|=1$. Let $R =\Gamma \setminus(\{0\}\cup I(\Gamma))$. For every positive integer $t$ and integer partition $\{m_i\}_{i=1}^t$ of $m-2$ with $m_i \geq 2$ for every $i\in[1,t]$, there is a subset partition $\{S_i\}_{i=1}^t$ of $R$ such that $|S_i| = m_i$ and $\sum_{s\in S_i}s = 0$ for every $i\in[1,t]$. \end{conj}

\section{Groups orthomorphism}
An \textit{orthomorphism} of a group $\Gamma$ is defined
as $\varphi\in$Bij$(\Gamma)$ (set of all bijections from $\Gamma$ to itself) that the mapping  $\theta\colon g \mapsto  g^{-1}\varphi(g)$ is
also bijective \cite{Hall}. 
Originally, orthomorphisms were introduced by Mann in 1942 as a tool for constructing mutually orthogonal Latin squares  \cite{Mann}. 

For finite Abelian groups, it was proved the following:
\begin{thm}[\cite{Hall,HP}]
A finite Abelian $\Gamma$ group has a complete
mapping if and only if $|I(\Gamma)|\neq1$.
\end{thm}
The most famous conjecture in the area is the Hall-Paige conjecture \cite{HP} which states that a group $\Gamma$ admits an orthomorphism if and
only if  $ \Gamma$ is a finite group and the Sylow 2-subgroups of $ \Gamma $ are either trivial or non-cyclic, then $ \Gamma $ has an orthomorphism.   Hall and Paige showed that a finite group with a nontrivial cyclic Sylow 2-subgroup does not admit complete mappings. The converse was established in 2009 by Wilcox 
\cite{Wilcox}, Evans \cite{Evans}, and Bray et al. \cite{Bray}.

A group $\Gamma$ of order $n$ is said to be $R$-\textit{sequenceable} if the nonidentity elements of the group can be listed in a sequence $g_1, g_2,\ldots, g_{n-1}$ such that  $g_1^{-1}g_2, g_2^{-1}g_3,\ldots, g_{n-1}^{-1}g_1$ are all distinct. This idea was introduced in  1974 by Ringel \cite{Ringel}, who used this concept in his solution of Heawood map coloring problem. Friedlander, Gordon and Tannenbaum generalized the notion of $R$-sequenceability by asking for which groups $\Gamma$ of order $n$, there exists an orthomorphism of $\Gamma$ fixing the identity and permuting the remaining elements as products of disjoint $k$-cycles for any $k$ dividing $n-1$ \cite{FrieGorTan}.
They stated the following conjecture:
\begin{conj}[\cite{FrieGorTan}] Let $\Gamma$ be an Abelian group of order $n$ scuh that $|I(\Gamma)|\neq 1$. Suppose for some integer $k \geq 2$ that $k$ divides $n - 1$. Then, there exists
an orthomorphism of $\Gamma$ that fixes the identity element, and permutes the remaining elements as products of disjoint cycles of length $k$.
\end{conj}

The above conjecture is still open. There are several partial results. Friedlander,
Gordon, and Tannenbaum confirmed their conjecture for groups of order at most 15, and Abelian $p-groups$ where $p \geq 3$ \cite{FrieGorTan}. Recently, the conjecture was confirmed for sufficiently large groups \cite{Muyesser}.

We will show the following:
\begin{thm}
Let $\Gamma$ be an Abelian group such that $|\Gamma|\cong 1\pmod{3}$, $|I(\Gamma)|\neq 1$. 
If $\Gamma\cong L\oplus H$, with $|L|=2^{\eta}\equiv 1\pmod 3$ for some positive integer $\eta$, and $|H|\equiv 1\pmod 6$, then there exists an orthomorphism of $\Gamma$ that fixes the identity element, and permutes the remaining elements as products of disjoint cycles of length $3$. 
\end{thm}
\textit{Proof.}  The group $\Gamma$ has $3$-ZSPP (see Theorem 2.2, \cite{CicrSuch2}), thus it can be written as:
$$\Gamma=\{0\}\cup\cup_{i=1}^{(|\Gamma|-1)/3}\{x_0^i,x_1^i,x_2^i\},$$
Where $x_0^i+x_1^i+x_2^i=0$ for any $i\in[1,(|\Gamma|-1)/3].$ Set now $\varphi(x_j^i)=-x_{j+2}^i$ for $j=0,1,2$, $i\in[1,(|\Gamma|-1)/3]$, where the subscripts are taken modulo $3$. Note that $\varphi(x_j^i)-x_j^i=-x_{j+2}^i-x_{j}^i=x_{j+1}^i$.~\qed

\section{Some applications in graph labeling}\label{sec:sa}
In this section, we explore applications of zero-sum sets in Abelian groups to problems involving magic-type labelings of graphs. Generally, such a labeling for a graph $G=(V,E)$ is defined as a mapping from either 
$V$, $E$, or their union $V\cup E$ to a set of labels, which is usually a set of integers or group elements. The weight of a graph element is typically calculated as the sum of labels of adjacent or incident elements, either of one type or both. When the weight of all elements is required to be identical, we refer to it as magic-type labeling, whereas if all weights must be distinct, it is called antimagic-type labeling. Perhaps the best-known problem in this field is the anti-magic conjecture by \cite{HR}, which posits that for every graph except 
$K_2$, there exists a bijective labeling of edges with integers $1, 2, \ldots, |E|$ such that each vertex has a unique weight. This conjecture remains unsolved.
\subsection{Antimagic labeling of trees}

In \cite{KLR}, Kaplan, Lev and Roditty considered the following generalization of the
concept of an antimagic graph.

Let $\Gamma$ be an Abelian group (not necessarily finite) an $A$ be a finite subset of $\Gamma^*=\Gamma\setminus\{0\}$ with $|A| = |E(G)|$. An $A$-labeling
of $G=(V,E)$ is a one-to-one mapping $f \colon E(G) \to A$. The weight of every vertex is calculated as the sum (taken in $\Gamma$) of the labels of incident edges.
 We shall say that $G$ is \textit{$A$-antimagic} if all the weights differ. In the case that $\Gamma$ is finite, we shall say that $G$ is $\Gamma^*$-antimagic if $G$ is $\Gamma^*$-antimagic.
They conjectured that a tree with $|\Gamma|$ vertices is $\Gamma$-antimagic if and only if $|I(G)|\neq 1$.

{A \textit{$k$-tree} $T$ is a rooted tree, where every vertex that is not a leaf has at least $k$ children.}
They showed the following
\begin{cor}[\cite{KLR}]
    Let $\Gamma$ be a finite abelian group with the $2$-ZSPP. Then
every 2-tree with $|\Gamma|$ vertices is $\Gamma^*$-antimagic
\end{cor}
Let $\{v_1,v_2,\ldots,v_t\}$ represent all the vertices of a 2-tree that are not leaves, with the corresponding numbers of children given by $\{r_1,r_2,\ldots,r_t\}$. Since $r_i\geq 2$, we can partition $\Gamma^*$ into zero-sum subsets $A_1,A_2,\ldots,A_t$, each with cardinalities 
$\{r_1,r_2,\ldots,r_t\}$, respectively. These subsets are then used to label the edges in the set $N_i=\{v_iw\colon w$ is a child of $v_i\}$ by elements of $A_i$. Given that every vertex in $T$ (except the root) has a unique parent, it follows that all vertex weights are distinct.

Kaplan, Lev, and Roditty proved that every $2$-tree has a $\zet_n$-antimagic labeling if $n$ is odd. Moreover, from this result follows that every $2$-tree of odd order is antimagic. Zeng proved that an Abelian finite group $\Gamma$ has $2$-ZSPP if and only if $|I(\Gamma)|\in\{0,3\}$ \cite{Zeng}. As a corollary he obtained that every $2$-tree has $\Gamma$-antimagic labeling if $|I(\Gamma)|\in\{0,3\}$. Theorem~\ref{SK2} implies that every $4$-tree of order $|\Gamma|$ is $\Gamma^*$-antimagic. 

 However, Kaplan et al. \cite{KLR} showed that, if $\Gamma$ has a unique involution, then any tree on $|\Gamma|$ vertices is not $\Gamma^*$-antimagic, the conjecture is not true. It is enough to consider a path $P_{2^m}$, $m>1$ and a group $\Gamma\cong(\zet_2)^m$. Since all non-zero elements in $\Gamma$ are involutions, there are not two of them that sum up to $0$. Note, that paths are not the only example of such trees, others are trees of order $n$ with maximum degree $n-2$ \cite{Cic24}. 

\subsection{Group irregular labeling of graphs}
Using the Pigeonhole Principle, it is easy to show that in any simple graph $G$, at least two vertices with the same degree exist. However, the situation changes if we consider an edge labeling $f:E(G)\rightarrow \{1, \ldots, k\}$ (where labels need not be distinct) and calculate the \textit{weighted degree} of each vertex $v$ as the sum of labels of all edges incident to $v$. The labeling 
$f$ is called \textit{irregular} if all vertex weighted degrees are unique, making it an antimagic-type labeling. The smallest value of $k$ that permits an irregular labeling is known as the \textit{irregularity strength of 
$G$}, denoted by $s(G)$. This problem was introduced by \cite{ref_ChaJacLehOelRuiSab1} and has since attracted considerable interest. Generally, it is known that $s(G) \leq n -1$ for any graph $G$ of order 
$n$ with no isolated edges and at most one isolated vertex, except for $K_3$, as shown by Aigner and Triesch \cite{AigTrie} and Nierhoff \cite{TigBouIrrStr}. Although the upper bound is tight for the family of star graphs, this result can be refined for graphs with sufficiently large minimum degree $\delta$. The best general upper bound was proven by Przybyło and Wei, who proved that for any $\varepsilon \in  (0, 0.25)$
there exist absolute constants $c_1, c_2$ such that for all graphs $G$ on $n$ vertices with minimum degree
$\delta\geq 1$ and without isolated edges, $s(G) \leq  n/\delta(1 + c_1/\delta^{\varepsilon}) + c_2$ \cite{PrzWei}.

Tuza began considering irregular labeling of a graph 
$G=(V,E)$ by the group $\Gamma=(\zet_2)^m$. In this labeling, edges received weights from elements of 
$\Gamma$ and the weighted degrees of vertices were computed using operations in $\Gamma$. Tuza denoted 
$m_M(G)$ as the smallest $m$ for which a $(\zet_2)^m$-irregular labeling of a graph $G$ exists, and showed that $m_M(G)\leq 3\log_2|V|$ \cite{Tuza}. This result was later improved by Aigner and Triesch in \cite{ref_AigTri94} to $m_M(G)\leq \left\lceil \log_2|V|\right\rceil+4$. These authors utilized the observation that this problem is related to subsets with zero-sum in the group $(\zet_2)^m$.

It is noteworthy that if there exists a $(\zet_2)^m$-irregular labeling of a graph $G$ with connected components $\{C_i\}_{i=1}^t$, then the sum of weighted degrees in each component equals zero (this follows because each edge's label is added twice to the weighted degrees). Aigner and Triesch conjectured $m_M(G)\leq \left\lceil \text{log}_2|V|\right\rceil+1$ for any graph $G$.  For large $|V|=n$, this was proven by Tuza using probabilistic methods in \cite{Tuza2}. It was subsequently proven for any graph by Caccetta and Jia \cite{CaccJia}, as well as independently by Egawa \cite{Egawa}, who focused on partitioning the group $\Gamma$ into subsets with zero-sum (we will put it into more details in the next section).

In \cite{ref_AnhCic1} other groups were considered. The authors introduced \textit{group irregularity strength} (denoted $s_g(G)$) of $G$ as the smallest $k$, such that for every Abelian group $\Gamma$ of order $k$ there exists a $\Gamma$-irregular labeling of $G$. Note that $s(G)\leq s_g(G)$ for every graph $G$.

Anholcer, Cichacz and Milani$\check{c}$ have shown the following theorem that describes the value of $s_g(G)$ for all connected graphs $G$ of order $n\geq 3$. 
\begin{thm}[\cite{ref_AnhCic1}]\label{Anh}
    Let $G$ be an arbitrary connected graph of order $n\geq 3$. Then
    $$s_g(G) = \left\{\begin{matrix}
        n+2 & \text{if } G\cong K_{1,3^{2q+1}-2} \text{ for some integer } q\geq 1\\ 
        n+1 & \text{if } n\equiv 2 (\text{mod} 4) \wedge G\ncong K_{1,3^{2q+1}-2} \text{ for any integer } q\geq 1 \\ 
        n & \text{otherwise} 
    \end{matrix}\right.$$
\end{thm}
It is easy to see that to distinguish all $n$ vertices in an arbitrary graph of that order we need at least $n$ different element of $\Gamma$.
Although the following lemma shows that an Abelian group of order $n$ is not always enough to have a $\Gamma$-irregular labeling of $G$.
\begin{lem}[\cite{ref_AnhCic1}]\label{2mod4}
    Let $G$ be a graph of order $n$.
    If $n\equiv 2 \pmod4$, then there is no $\Gamma$-irregular labeling $G$ for any Abelian group $\Gamma$ of order $n$.
\end{lem}
So far the best upper bound for any graph was given by Anholcer, Cichacz, and Przybyło in \cite{LinearBoundNo0}
\begin{cor}[\cite{LinearBoundNo0}]\label{linear}
    Let $G$ be an arbitrary graph of order $n$ having no component of the order less than 3. 
    Then  $s_g(G) \leq 2n.$
\end{cor}

In \cite{ref_AnhCic} used $2$-ZSPP of groups of odd order for bounding the group irregularity strength of disconnected graphs without a star as a connected component. Roughly speaking, the authors divide every connected component into $2$-subsets and $3$-subsets of vertices and partition the set of non-zero elements of the corresponding group into the same number of zero-sum $2$-subsets and $3$-subsets, and later use the method of augmented paths to do the labeling.
Namely, for any two given vertices $v_1$ and $v_2$ from the same connected component of a graph $G$ exist walks from $v_1$ to $v_2$. While labeling edges of a graph $G$ they started with $0$ on all of them.
Next, they choose $v_1$ and $v_2$ and modify all of the edges of a chosen walk from  $v_1$ to $v_2$ by adding some element of an Abelian group $\Gamma$.
Now they add some element $a\in \Gamma$ to all the labels of the edges in an odd position on the walk, starting from $v_1$ and $-a$ to the rest of the labels of edges. 
This modification of labels increases the weighted degree of $v_1$ by $a$ and increases the weighted degree of $v_2$ by $a$ or $-a$. We should also note that in both cases the operation does not modify the weighted degree of any other vertices of the walk.

In \cite{CicKrup} the authors use the above method and Theorem~\ref{Tannenbaum1} to improve Theorem~\ref{linear} for graphs without small stars.
\begin{thm}[\cite{CicKrup}]\label{bez_K13}
    Let $G$ be a graph of order $n$ having neither a component of the order less than 3 nor a $K_{1,1+2u}$  component for $u\in\{1,2\}$.
    Then
     $$\begin{matrix}
       s_g(G) = n & \text{if} & n\equiv 1 \pmod2\\ 
       s_g(G) \leq n+3 & \text{if} & n\equiv 0 \pmod2.
    \end{matrix}$$
\end{thm}

Moreover it was shown that for a graph $G$  of order $n$ having neither a component of the order less than 3 nor a $K_{1,1+2u}$  component for $u\in\{1,2\}$ there exists a $\Gamma$-irregular labeling in any group of odd order $t\geq n+3$.

\subsection{Group irregular labeling of directed graphs}
Let  $\overrightarrow{G}$  be a directed  graph of order $n$. If  there exists a
 mapping $\psi$ from  $E(\overrightarrow{G})$  to an Abelian group $\Gamma$
such
that if we define a mapping $\varphi_{\psi}$ from $V(\overrightarrow{G})$ 
to $\Gamma$
by
$$\varphi_{\psi}(x)=\sum_{y\in N^+(x)}\psi(xy)-\sum_{y\in N^-(x)}\psi(yx),\;\;\;(x\in V(\overrightarrow{G})),$$
then $\varphi_{\psi}$
is injective, then such a labeling $\psi$ is called $\Gamma$-\textit{irregular}. In this situation, we say that $\overrightarrow{G}$ is \textit{realizable} in
$\Gamma$.

Since for $\Gamma=(\zet_2)^m$, it does not matter whether we consider $\Gamma$ irregular labeling of a simple graph $G$ or a directed graph $\overrightarrow{G}$, this problem for $(\zet_p)^m$
  (where $p$ is an odd prime number) was defined by Fukuchi in \cite{Fukuchi} as a generalization of the above-mentioned problem of Tuza \cite{Tuza}. Analogously to the results in \cite{CaccJia, Egawa}, Fukuchi showed that there exists a  $(\zet_p)^m$-irregular labeling of the directed graph $\overrightarrow{G}$   with weakly connected components $\{\overrightarrow{G}_i\}_{i=1}^t$
  if and only if there are pairwise disjoint subsets 
 $\{S_i\}_{i=1}^t$   in  $(\zet_p)^m$   such that 
 $|S_i|=|V(\overrightarrow{G}_i)|$ and $\sum_{s\in S_i}s=0$ for $i\in[1,t]$. This result was then extended to any Abelian group $\Gamma$ by Cichacz and Tuza  \cite{CicTuz}.
 
 Fukuchi demonstrated that 
$(\zet_p)^m$ can be partitioned into zero-sum subsets $\{S_i\}_{i=1}^t$, thereby achieving the result stating that any directed graph $\overrightarrow{G}$
  of order $n$ is realizable in any group $(\zet_p)^m$
  provided $n\leq p^m$. 

Based on the results of Tannenbaum \cite{Tannenbaum1} and Zeng \cite{Zeng}, we immediately conclude that if $|I(\Gamma)|\in\{0,3\}$ and $|\Gamma|\geq |V(\overrightarrow{G})|$, then any 
 $\overrightarrow{G}$
  without components of order less than $3$ is realizable in  $\Gamma$.

Moreover, using the result from \cite{Cicir,CicrSuch,CicrSuch2} we obtain the following.
\begin{thm}[\cite{Cicir,CicrSuch,CicrSuch2}] Any digraph\/ $\overrightarrow{G}$ of order\/ $n$ with no weakly connected components
 of order less than\/ $4$ has a $\Gamma$-irregular labeling for every\/ $\Gamma$ such that\/
   $|\Gamma|\geq n+6$.\end{thm}
Cichacz and Tuza proved that if  $n$ is large enough
with respect to an arbitrarily fixed $\varepsilon > 0$ then $\overrightarrow{G}$
has a $\Gamma$-irregular labeling for any $\Gamma$ such that $|\Gamma|>(1+\varepsilon)n$ \cite{CicTuz}. 
Therefore Cichacz and Suchan stated the conjecture:

\begin{conj}[\cite{CicrSuch2}]\label{hipotezaK}There exists a constant $K$ such that any digraph\/ $\overrightarrow{G}$ of order\/ $n$ with no weakly connected components
 of order less than\/ $3$ has a $\Gamma$-irregular labeling for every\/ $\Gamma$ such that $|\Gamma|\geq n+K$.
\end{conj}

\subsection{Group distance magic labeling}\label{subsec:gdm}

A $\Gamma$\emph{-distance magic labeling} of a graph $G = (V, E)$ with $|V| = n$ is a bijection $\ell$ between $V$ and  an Abelian group $\Gamma$ of order $n$
such that the weight $w(x) =\sum_{y\in N(x)}\ell(y)$ of every vertex $x \in V$ is equal to the same element $\mu\in \Gamma$, called the \emph{magic
constant}.

Notice that the constant sum partitions of a group $\Gamma$ lead to complete multipartite $\Gamma$-distance magic labeled graphs. For
instance, the partition $\{0\}$, $\{1, 2, 4\}$, $\{3, 5, 6\}$ of the group $\zet_7$ with constant sum $0$ leads to a $\zet_7$-distance magic labeling
of the complete tripartite graph $K_{1,3,3}$ (see \cite{CicZ}). Indeed, suppose we have a constant sum partition of $\Gamma$: $\{S_i\}_{i=1}^t$ of $\Gamma$ with $|S_i| = m_i$ and $\sum_{s\in S_i}s = \nu$ for every {$i\in[1,t]$} and some $\nu \in \Gamma$. Let $G$ be a complete $t$-partite graph with the color classes $\{A_i\}_{i=1}^t$, where $|A_i|=m_i$ for every  {$i\in[1 ,t]$}. Let us label the vertices of $A_i$ with distinct elements of $S_i$ for every {$i\in[1 ,t]$}, and compute the weight of every vertex as the sum (in $\Gamma$) of the labels of its neighbors. It is easy to see that all the weights are equal and thus we have a $\Gamma$-distance magic labeling. On the other hand, suppose $G=K_{m_1,\ldots,m_t}$ is a complete $t$-partite graph of order $m$ with the color classes $\{A_i\}_{i=1}^t$ that is $\Gamma$-distance magic  {with a labeling $\ell$}. So $\sum_{i=1,i\neq j}^t\sum_{x\in A_i}\ell(x) = \mu$ for every {$j\in[1 , t]$}, which implies that $\sum_{x\in A_j}\ell(x) = \nu$ for every  {$j\in[1 , t]$}, and some $\nu\in\Gamma$.

Let $G=K_{m_1,\ldots,m_t}$ be a complete $t$-partite graph of order $m$. Let now $1\leq m_1\leq m_2 \leq\ldots\leq m_t$. Using some constant-sum partition properties of every finite Abelian group $\Gamma$ of order $m$, it was shown that:
\begin{enumerate}
\item  For $t=2$, if  $m_1 + m_2\not\equiv 2\pmod 4$ then the graph $G$ is  $\Gamma$-distance magic \cite{CichaczOM}.
	\item  For $t=3$, if ($m_2>1$ and  $m_1 + m_2+m_3\neq 2^p$ for any positive integer $p$) or  ($m_1\neq 2$ and $m_2>2$), then the graph $G$ admits a $\Gamma$-distance magic labeling \cite{Cic3}.
	\item If $m_1=m_2=\ldots=m_t>2$ and $|I(\Gamma)|\neq 1$, then the graph $G$ admits a $\Gamma$-distance magic labeling \cite{CicZ}.
	\item If $m_1\geq 3$ and $m_t\geq \frac{1}{2}(m +\sqrt{2m+1}) -1$, and $|I(\Gamma)|\neq 1$, then\/ $G$ admits a $\Gamma$-distance magic labeling \cite{CicTuz}.	
	\item  If $m_1\geq 4$, $m$ is large enough, and $|I(\Gamma)|\neq 1$, then\/ $G$ admits a $\Gamma$-distance magic labeling \cite{CicTuz},
	\item If $m_2\geq 2$, $t$ is odd, and $\Gamma\cong\zet_n$, then\/ $G$ admits a $\Gamma$-distance magic labeling \cite{FREYBERG2019}.
 	\item If $m_1\geq 3$ and $m=2^n$ for some $n$, and $\Gamma\not\cong \zet_{2^n}$, then\/ $G$ admits a $\Gamma$-distance magic labeling \cite{CicrSuch}.
\end{enumerate}
Moreover, Theorem~\ref{SK2} implies that if $m_1\geq 4$ and $|I(\Gamma)|\neq 1$, then\/ $G$ admits a $\Gamma$-distance magic labeling. 


\section{Statements and Declarations}
This work was partially  supported by program ''Excellence initiative – research university'' for the AGH University.

\bibliographystyle{abbrv}
\bibliography{bibliografia1}

\end{document}